\documentclass[11pt,fleqn]{article}

\usepackage{amsmath}
\usepackage{amssymb}
\usepackage{times}
\usepackage{mathptm}
\usepackage{ifthen}
\usepackage{graphics}
\usepackage{color}
\usepackage{subfigure}
\usepackage{epsfig}
\usepackage[pdftex,colorlinks,linkcolor={blue},citecolor={blue}]{hyperref}
\usepackage[margin=1in]{geometry}
\usepackage[T1]{fontenc}
\usepackage[authoryear,round]{natbib}
\usepackage{titlesec}
\usepackage{appendix}
\usepackage{soul}

\newcommand{\1}{\textrm{\textbf{1}}}

\newcommand{\ta}{\tau}

\newcommand{\tb}{\bar{\tau}}
\newcommand{\eib}{\bar{e} (i)}
\newcommand{\cib}{\bar{c} (i)}
\newcommand{\cjb}{\bar{c} (j)}
\newcommand{\bib}{\bar{b} (i)}
\newcommand{\bjb}{\bar{b} (j)}
\newcommand{\yib}{\bar{y} (i)}
\newcommand{\yjb}{\bar{y} (j)}
\newcommand{\eb}{\bar{e}}
\newcommand{\cb}{\bar{c}}
\newcommand{\yb}{\bar{y}}

\renewcommand{\1}{\textrm{\textbf{1}}}

\newtheorem{theo}{Theorem}
\newtheorem{lemma}{Lemma}
\newtheorem{corollary}{Corollary}

\newtheorem{fact}{Fact}
\newtheorem{remark}{Remark}
\long\def\comment#1{}

\begin{document}

\title{Analysis of Fixed-Time Control\thanks{This research was supported in part by NSF SBIR Award 1329477 and by the California Department of Transportation under the Connected Corridors program.  The authors are  grateful to Professors Hong Lo, Ravi Mazumdar and Jean Walrand for their very helpful suggestions.} }
\author{Ajith Muralidharan, Sensys Networks, Inc, Berkeley, CA 94710 \\
amuralidharan@sensysnetworks.com\\ Ramtin Pedarsani,
University of California, Berkeley, CA 94720-1770\\
{ramtin@berkeley.edu}\\
Pravin Varaiya, University of California, Berkeley, CA 94720-1770\\ Corresponding Author, varaiya@berkeley.edu, 1-510-642-5270
 }
\date{}
\maketitle
\begin{abstract}

The paper presents an analysis of the traffic dynamics in a network of signalized intersections. The intersections are regulated by fixed-time (FT) controls, all with the same cycle length  or period, $T$. The network is modeled as a queuing network.    Vehicles arrive from  outside the network at entry links in a deterministic periodic stream, also with period $T$. They take a fixed time to travel along each link,
and at the end of the link they join a queue.  There is a separate queue at each link for each movement or phase.  Vehicles  make turns at intersections in fixed proportions,  and eventually leave the network, that is, a fraction $r(i,j)$ of vehicles that leave queue $i$  go to queue $j$ and the fraction $[1 - \sum_j r(i,j)]$ leave the network. 
The storage capacity of the queues is infinite, so there is no spill back. The main contribution of the paper is to show that if the signal controls accommodate the demands then, starting in any initial condition, the network state converges to a unique periodic orbit. Thus, the effect of initial conditions disappears. More precisely, the state of the network at time $t$ is the vector $x  (t)$ of all queue lengths, together with the position of vehicles traveling along the links.  Suppose that the network is stable, that is,  $x(t)$ is bounded. Then \\
\hspace*{0.2in}(1) there exists a unique periodic trajectory $x^*$, with period $T$;\\
\hspace*{0.2in}(2) every trajectory converges to this periodic trajectory; \\
 \hspace*{0.2in}(3) if vehicles do not follow loops, the convergence occurs in finite time. \\
The periodic trajectory determines the performance of the entire network.


\end{abstract}
\section{Introduction}\label{intro}
Traffic in an urban network  is determined by intersection signal control and the pattern of demand.    The movement of vehicles is often modeled as  a queuing network as, for example, in \cite{papagIEEE,rhodes}.  Roughly speaking,  a vehicle arrives from outside the network at an entry link; travels along a link at a fixed speed; at the end of the link it arrives at an intersection and joins a queue of vehicles for the next link in its path; the queue is served at a specified saturation flow rate when that movement is actuated by the signal; eventually the vehicle leaves the network.

In the U.S. 90 percent of traffic signals follow fixed time (FT) controls, which operate the signal in a fixed periodic cycle, independent of the traffic state (\cite{stm}).   Despite its practical importance,  little attention has been paid to understanding how traffic behaves under under FT control.  Published work has studied queues at a single, isolated intersection, as in \cite{miller63}.  The steady state optimal control of single  intersections is studied in \cite{Improta,Jhaddad,GazisBook}. The latter work derives the optimal control settings required to minimize different objectives including queuing delays, but does not address the effect of initial conditions on solution trajectories or their convergence. \cite{GazisMultiIntersectionsOversaturated,GazisBook} analyze oversaturated intersections and \cite{MPtrc} inroduces an adaptive control for undersaturated networks.  But neither work  analyzes the behavior of solution trajectories. Signal timing tools used by  traffic engineers often employ empirical models  (\cite{Webster,HCM2010}) in combination with simulations,  assuming steady state conditions.  But the absence of theory establishing  convergence to a unique steady state calls into question whether the traffic flows  achieve the performance for which these signals are tuned.

We  analyze  vehicle movement under two assumptions: first, all the signals have a fixed time (FT) control with the same cycle time or period $T$; second, vehicles from outside enter the network in periodic streams with the same period.   Periodic demands include constant demands, which is the assumption  in commercial packages used to design FT controls.  Also, if there are intersections with FT controls with different cycles $T_1, \cdots, T_k$,  they are  all periodic with the same period $T = \text{lcm} \{T_1, \cdots, T_k\}$.

The state of the signalized network at any time $t$ consists of $x(t)$, the vector of all queue lengths, together with the position of all vehicles that are traveling along a link but have not yet reached a queue.  A  queue increases when vehicles arrive and decreases when the control serves that queue.
We treat time as continuous and vehicles as a fluid instead of as discrete entities. As a result the evolution of the network is described by a delay-differential equation, in which the delay comes from the travel time of a vehicle as it moves from one queue to the next.  In an actual transportation network, the arrival and service processes are stochastic. However, an exact analysis of queue-length processes in a stochastic queueing network is very difficult, if not impossible; except for very simple examples such as an isolated intersection, the
underlying Markov chain of the system is intractable. Therefore, we consider deterministic arrival and service processes in this paper.

From a traffic theory viewpoint, our main contribution is to show that there is a unique periodic trajectory $x^*(t)$ of the queue length vector to which every trajectory $x(t)$ converges; moreover, in case individual vehicles do not circulate in loops,  the convergence is in finite time.  The periodic orbit of course determines every possible performance measure, such as
delay, travel time, amount of wasted green, and signal progression quality, see \cite{bullock-performance}.  An outstanding open problem is to calculate this periodic orbit without simulation.  If  this can be done, one would have a computational procedure to design the FT control for a network that optimizes any performance measure.

The results have some independent mathematical interest.  The delay-differential equation is not Lipschitz, and  existence and uniqueness of a solution is established using the reflection map of queuing theory (\cite{harrison-reiman,reflect-whitt}).   The differential equation is periodic (with period $T$), and the existence of a periodic orbit is proved using the Poincare map.  The global stability of this periodic orbit depends on  a monotonicity property reminiscent of that in freeway models (\cite{gomes07}).

The rest of the paper is organized as follows.  \S \ref{sec-single} presents the main results for a single queue.  \S \ref{sec-network} describes the basic results for the
network model.  Results for the case of periodic demand and FT control are presented in \S \ref{sec-periodic}.  The main conclusion and
some open questions are summarized in \S \ref{sec-conc}.

\section{Single queue without routing} \label{sec-single}
Time is continuous, $t \ge 0$.  The length or size of a single queue $x(t), t\ge 0$, evolves as
\begin{equation}
\dot{x}(t)= e(t) - b(t),  \label{1}
\end{equation}
with  arrivals $e(t) \ge 0$,  departures $b(t)$, $t \ge 0$, and initial queue $x(0) = x_0 \ge 0$.  Arrivals $e(t)$ are exogenously specified. There is a specified saturation flow  or service rate $c(t) \ge 0, t \ge 0$, so departures are  given by
\begin{equation} \label{2}
b(t) = \left \{
\begin{array}{ll }
c(t), & \mbox{ if } x(t) > 0, \\
\in [0, c(t)], & \mbox{ if } x(t) = 0,\\
0, & \mbox { if } x(t) < 0.
\end{array}
\right .
\end{equation}
Express the departure process as
\begin{eqnarray}
b(t) = c(t) - y(t), \; t \ge 0, \label{3}
\end{eqnarray}
so $y(t)$ is the rate at which  service  is unused.  From \eqref{2},
\begin{eqnarray*}
y(t) \ge 0, \mbox{ and } x(t) y(t) \equiv 0.
\end{eqnarray*}
Rewrite \eqref{1} as
\[\dot{x}(t) = [e(t) - c(t)] + y (t) ,\]
or in functional form as
\begin{equation}
x = u + v ,\label{4}
\end{equation}
in which
\begin{equation}
u(t) = x_0 + \int_0^t [e(s)-c(s)] ds \; \mbox{ and }v(t) = \int_0^t y(s) ds. \label{5}
\end{equation}
Then $x, v$ satisfy
\begin{equation}
x(t) \ge 0, \; v(t) \ge 0, \; v(0) = 0, \; y(t) = \dot{v}(t) \ge 0, \mbox{ and } x(t) \dot{v}(t) \equiv 0 . \label{6}
\end{equation}
Observe that $v(t)$ has an interpretation as the cumulative unused service (wasted green).


Fact \ref{F1} and Theorem \ref{T1} are immediate consequences of \cite[Theorem 1]{harrison-reiman}.

\begin{fact}\label{F1}
Suppose $x, b$ satisfy \eqref{1}-\eqref{2}.  Define $u, v$ by \eqref{3} and \eqref{5}.  Then \eqref{4} and \eqref{6} hold.
Conversely, suppose $u, v, x$ satisfy \eqref{4}-\eqref{6}.  Define $b (t)$ by \eqref{3}.  Then \eqref{1}-\eqref{2} hold.
\end{fact}
\textbf{Proof} \ Suppose $x, b$ satisfy \eqref{1}-\eqref{2} and define $u, v$ by \eqref{3} and \eqref{5}.   Then \eqref{5}  implies \eqref{4} and \eqref{2} implies \eqref{6}.
Conversely suppose \eqref{4}-\eqref{6} hold.  Define $b(t)$ by \eqref{3}.  Then, by \eqref{4},
\[\dot{x}(t) = \dot{u}(t) + y(t) = e(t) - c(t) + y(t) = e(t) - b(t),\]
so \eqref{1} holds.  Moreover, by  \eqref{6}, $x(t) > 0$ implies $y(t) = c(t) - b(t) = 0$, so $b(t) = c(t)$.  If $x(t) = 0$, $t \in (t_1, t_2)$, $\dot{x}(t) = 0, 0 \le y(t) = c(t) - e(t) \le c(t)$
and $b(t) = c(t) - y(t) = e(t) \le c(t)$.  This proves \eqref{2}. \hfill $\Box$

\begin{theo}\label{T1}
Fix continuous function $u$ with $u(0) =x_0 \ge 0$.  There exist unique continuous functions $x, v$ satisfying \eqref{4}-\eqref{6}. The functions  $v = \psi(u) $ and $x = \phi(u)$ are continuous (in sup norm) and given by
\begin{eqnarray}
v(t) &=&   \sup \{[u(s)]^- ~|~ 0 \le s \le t\}, \mbox{ with } [z]^- = \max\{ -z, 0\} ,\label{7} \\
x(t) &=& u(t) + \sup \{[u(s)]^- ~|~ 0 \le s \le t\} .\label{8}
\end{eqnarray}
\end{theo}
\textbf{Proof} \
Fix continuous  $u$ with $u(0) = x_0 \ge 0$. \\ 
\textit{Existence} \
 Define
\begin{equation}
v(t) = \sup_{0\le s \le t} [u (s)]^-  \mbox{ and } x(t) = u(t) + v(t) .\label{A1}
\end{equation}
Then \eqref{4} holds.  $v(0) = [x_0]^- = 0$ and $v$ is increasing.  Further, $x(t) \ge 0$, since
$x(t) = u(t) +\sup_{0\le s \le t} [u (s)]^- \ge 0$ if $u(t) \ge 0$, and $x(t) \ge u(t) + [-u(t)] \ge 0$ if $u(t) < 0$.
Suppose $\dot{v}(t) > 0$.  Then $v(t) = \sup_{0\le s \le t} [u (s)]^- = [u(t)]^- = -u(t)$,  so $v(t) + u(t) = 0 = x(t)$.  This proves \eqref{6}.

\textit{Uniqueness} \
Consider any solution
\begin{equation}
q= u+w, \; q \ge 0, \; w\ge 0, \; w(0) = 0, \; \dot{w} \ge 0, \mbox{ and } q \dot{w} = 0.  \label{A2}
\end{equation}
Since $w \ge 0$ and increasing,
\[w(t) = [q(t)-u(t)] = \sup_{0\le s \le t} [q(s) - u(s)]\ge \sup_{0 \le s \le t} [u(s)]^- = v(t).\]
If $w(t) > v(t)$, then there is $t_0 < t$ with $w(t_0) > v(t_0)$ and $\dot{w}(t_0) > 0$.  But then
\[q(t_0) = w(t_0) + u(t_0) > v(t_0) + u(t_0) = x(t_0) \ge 0.\]
So $q(t_0) > 0$ and $\dot{w}(t_0) > 0$, which contradicts \eqref{A2}.  So $w \equiv v$.

\textit{Continuity}\
Suppose $v = \psi (u)$ and $v' = \psi (u')$, and $\sup_{0\le s \le t} |u(s)-u'(s)| < \epsilon$.  Then
\[\sup _{0\le s \le t} [u(s)]^- \le \sup _{0\le s \le t} [u'(s)]^- + \epsilon,\]
and so $|v(t)-v'(t)| < \epsilon$.   Hence $\psi$ is continuous and so is $\phi$.
 \hfill $\Box$

Corollary \ref{C1} describes a useful monotonicity property.
\begin{corollary} (Monotonicity) \label{C1}
(a) If $u \le u'$ (pointwise), then $v = \psi(u) \ge v' = \psi (u')$ .\\
(b) If $x_0 \le x'_0 $, $e (t) \le e'(t)$ and $c(t) = c'(t)$, for all $t$,  $x (t)\le x' (t)$ and $b (t)\le b'(t)$, for all $t$.
\end{corollary}
\textbf{Proof} \ (a) If $u \le u'$, $[u(s)]^- \ge [u'(s)]^-$, and so $v \ge v'$, \\
(b) It is enough to show that $x(0) = x'(0)$ implies $x(t) \le x'(t)$ for small $t$.  If $x(0) = x'(0) > 0$ then, from \eqref{1}-\eqref{2},
\[\dot{x}(t) = e(t) - c(t) \le e'(t) - c'(t) = \dot{x}'(t),\]
so $x(t) \le x'(t)$ for small $t$.

If $x(0) = x'(0) = 0$, recall that
\[u(t) = \int_0^t [e(s) - c(s)]ds, \quad u'(t) = \int _0^t [e'(s) - c'(s)]ds.\]
Then $\dot{u}(t) \le \dot{u}'(t)$ and
\[x(t) = u(t) + \sup _{0\le s \le t}[u(s)]^- = u(t) + [u(s)]^- \mbox{ for some } s .\]
If $u(s) < 0$,
\begin{eqnarray*}
x(t)
&=& \int_0^t \dot{u}(\tau) d\tau - \int_0^s \dot{u}(\tau) d\tau= \int_s^t \dot{u}(\tau) d\tau \\
&\le & \int_s^t \dot{u'}(\tau) d\tau = \int_0^t \dot{u'}(\tau) d\tau - \int_0^s \dot{u'}(\tau) d\tau\\
&\le & u'(t) + \sup _{0\le s \le t}[u'(s)]^- = x'(t).
\end{eqnarray*}
If $u(s) \ge 0$ for $s \le t$, then $x(t) = u(t) \le u'(t) = x'(t)$.

Lastly, if $x'(t) > 0$ then from \eqref{2} $b'(t) = c'(t) = c(t) \ge b(t)$ and if $x'(t) = x(t) = 0$, $t \in (t_1, t_2)$,  $\dot{x}(t) = \dot{x}'(t) = 0$, so from \eqref{2},
\[b(t) = e(t) \le e'(t) = b'(t).\]
 \hfill $\Box$
\begin{corollary} \label{C2}
Suppose $u(0) = 0$ and $u(t) \rightarrow -\infty$.  Consider two solutions $x, x'$ of \eqref{4} with
\[x(t) = x_0 + u(t) + v(t), \quad x'(t) = x'_0 + u(t) + v'(t).\]
Then there exists $T^*$ such that $x(t) = x'(t), \ t > T^*$.
\end{corollary}
\textbf{Proof}
Suppose $x_0 < x'_0$. By Corollary \ref{C1}, $x(t) \le x'(t)$ for all $t$.  Suppose $x(t) < x'(t)$ for all $t$.  Then $x'(t) > 0$ for all $t$ and so from \eqref{6} $v'(t) \equiv 0$,
but then $x'(t) = x'_0 + u(t) \rightarrow - \infty$, which contradicts $x'(t) \ge 0$.  So there exists $T^*$ such that $x(T^*) = x'(T^*)$, and then $x(t) = x'(t)$ for $t \ge T^*$. \hfill $\Box$

\begin{remark} \label{R1}
The condition
$u(t)\rightarrow -\infty$ means that
\[\int_0^t [c(s) - e(s)] ds \rightarrow \infty,\]
which holds if on average the service rate exceeds the  arrival rate by some $\epsilon > 0$.  In turn, Corollary \ref{C2} says that the effect of initial condition $x_0$
disappears after a finite time.
\end{remark}

\begin{theo} \label{T2}
 Suppose $u(0) = 0$, $\dot{u}(t)$ is periodic with period $T$, and $\int _0^T \dot{u}(t)dt < 0$.  Then there is a unique periodic trajectory $z$ with period $T$ such that
 \begin{eqnarray*}
 z(t) = z_0 + u(t) + v(t).
 \end{eqnarray*}
 Every solution of
 \begin{equation}\label{11}
 x(t) = x_0 + u(t) + v(t),
 \end{equation}
 coincides with $z(t)$ after some finite time.  There exists $t_0 \in [0,T]$ such that $z(t_0) = 0$, i.e. the queue will be cleared in each period.
 \end{theo}
\textbf{Proof} \
Consider the Poincare map $F: x(0) \mapsto x(T)$ and the iterates $x(nT) = F^n (x(0))$.  Take $x(0) = 0$. Clearly, $x(T) \geq x(0) = 0$. Hence, by monotonicity,
$$
x(0) \leq x(T) \leq x(2T) \leq \cdots .
$$
Since $\int_0^T \dot{u}(t) dt < 0$, the queue length is bounded.  Since the sequence is bounded and increasing, it converges. The uniqueness and convergence in finite time are immediate results of Corollary \ref{C2}. Finally, we show that the queue will be cleared in each period. Let $z$ be the unique periodic trajectory. If $z(t)= z_0 + u(t) + v(t) > 0$ for all $t$, then $v(t) \equiv 0$  and since $u(t) \rightarrow - \infty$, $z(t) \rightarrow - \infty$, but $z(t) \ge 0$.   So $z(t_0) = 0$ for some $t_0 $.
\hfill $\Box$

\textbf{Example 1} \
\begin{figure}[h!]
\centering
\includegraphics[width = 4.5in]{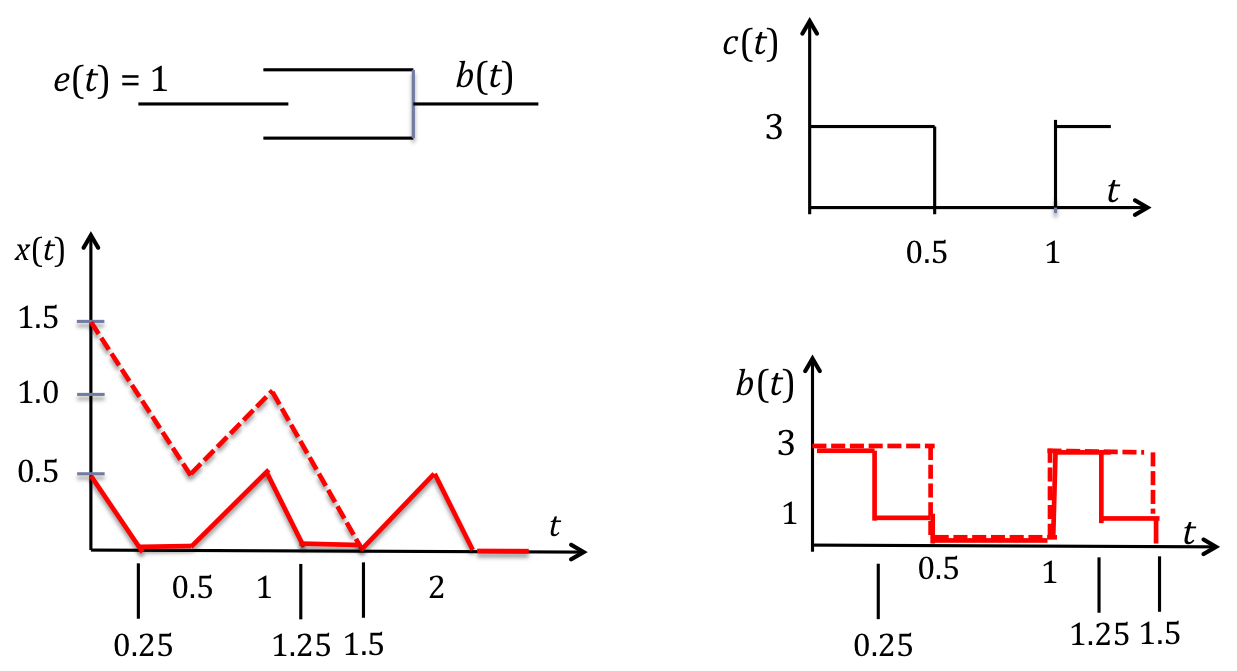}
\caption{Example1: constant arrival rate, $e(t) = 1$; periodic service rate $c(t)$ with period 1; two solutions converge at $t = 1.5$.}
\label{fig1}
\end{figure}
Consider a single link with constant arrival rate $e(t) = 1$.  The cycle time is 1, and the periodic saturation rate is $c(t) = 3$ for $0 \leq t \leq 0.5$, and $c(t) = 0$ for $0.5 < t < 1$. Thus within each cycle the signal is green for time 0.5, and red for time 0.5.   During green, 3 vehicles can depart per unit time.  Figure \ref{fig1} shows the unique periodic orbit (solid line) starting at $x(0) = 0.5$ and another trajectory
(dashed line) starting at $x(0) = 1.5$, which coincides with the periodic orbit after $t = 1.5$.  Also shown are the two associated departure
processes $b(t)$, which also coincide after $t=1.5$.

\section{Network of queues}\label{sec-network}
\begin{figure}[h!]
\centering
\includegraphics[width=4.0in]{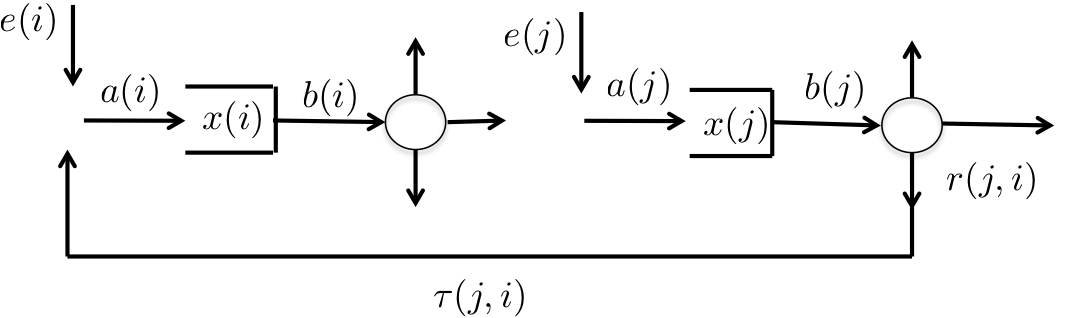}
\caption{Evolution of $x(i)(t) $.}
\label{fig-network}
\end{figure}
Figure \ref{fig-network} will help establish notation for  the network.  A fraction $r(j,i)$ of vehicles leaving queue $j$  will travel along link $(j,i)$ and join queue $i$ after time $\ta (j,i)$, t  Vehicles join queue $i$ at rate $a(i)$ either from outside the network at rate $e(i)$ or  after being
routed from another queue.  It is assumed that each queue has infinite storage capacity.

Hence the queuing equations for the network in Figure \ref{fig-network} are
\begin{eqnarray*}
\dot{x}(i)(t)& =& a(i)(t) - b(i)(t) , \ t \ge 0, \\
a(i)(t) &=& e(i)(t) + \sum_j b(j)(t-\ta(j,i)) r(j,i) ,\\
b(i)(t) &= &
\left \{
\begin{array}{ll}
c(i) (t) , & \mbox{ if } x(i)(t) > 0, \\
\in [0, c(i) (t)], & \mbox{ if } x(i)(t) = 0,\\
0, & \mbox { if } x(i)(t) < 0.
\end{array}
\right .
\end{eqnarray*}
Above $c(i)(t)$ is the saturation flow or service rate at which queue $i$ is served.

Express the departure process as
\begin{equation} \label{12}
b(i)(t) = c(i)(t) - y(i)(t), t \ge 0.
\end{equation}
Then $y(i)(t) \in [0, c(i)(t)]$ is the unused service rate,
\begin{equation} \label{13}
y(i)(t) \ge 0, \mbox{ and } x(i)(t) y(i)(t) \equiv 0.
\end{equation}
Rewrite the system equations as
\begin{equation}
\dot{x}(i)(t)= e(i)(t) - c(i)(t) + \sum_j b(j)(t-\ta (j,i)) r(j,i)  + y(i)(t)  , \label{14}
\end{equation}
or in functional form as
\begin{eqnarray}
x(i) = u(i) + v(i), \label{15}
\end{eqnarray}
in which
\begin{eqnarray}
u(i)(t) &=& x(i)(0) + \int_0^t \{e(i)(s) - c(i)(s) + \sum_j b(j) (s-\ta (j,i)) r(j,i)\} ds, \label{16}\\
b(j) (s) &=& c(j)(s)- y(j)(s) \label{17}\\
v(i)(t) &=& \int_0^t y(i)(s) ds.\label{18}
\end{eqnarray}
Fix the external arrivals and saturation flows, $\{e(i)(t), c(i)(t), t \ge 0\}$.  Suppose $0 < \ta = \min\{\ta (i,j)\} \le \tb = \max \{\ta (i,j)\}$.  Fix  initial conditions
\begin{equation} \label{19}
x(i) (0)\ge 0 \mbox { and }\{ b(i)(s) \ge 0, s \in [-\tb, 0]\}.
\end{equation}
 This determines $u(i)(s), s \in [0, \ta]$.  By Theorem \ref{T1}, there exist unique $\{x(i)(t), v(i)(t)\}$ satisfying \eqref{13}-\eqref{14} for $t \in [0, \ta]$.  In turn this fixes
new initial conditions at time $\ta$ similar to \eqref{19}:
\[ \{x(i)(\ta)\ge 0\} \mbox{ and } \{b(i)(s) = c(i)(s) - y(i)(s)\ge 0, s \in [-\tb + \ta, \ta]\}, \]
so that again by Theorem \ref{T1} the solution can be extended to $[\ta, 2 \ta]$.  Proceeding step wise in this way leads to the next result.
\begin{theo}\label{T2}
Fix arrivals $\{e(i)(t)\}$, saturation rates $\{c(i)(t)\}$, and routing ratios  $\{r(i,j)\}$.  Fix initial conditions $\{x(i)(0) \ge 0\}, \{b(i)(s), s \in [-\tb, 0]\}$.
Then there are unique  functions $\{x(i)(t), b(i)(t), v(i)(t)$, $ t \ge 0\}$ satisfying \eqref{12}-\eqref{14}.
\end{theo}

\begin{corollary} (Monotonicity)\label{C3}
Suppose $x(i)(0) \le x'(i)(0), 0 \le b(i)(s) \le b'(i)(s), s \in [-\tb, 0]$, $e(i)(t) \le e'(i)(t)$, $c(i)(t) =
\begin{scriptsize}
{\footnotesize •\begin{small}
{\normalsize •}
\end{small}}
\end{scriptsize} c'(i)(t)$, $t \ge 0$. Then
$x(i)(t) \le x'(i)(t)$, $b(i)(t) \le b'(i)(t)$, $u(i)(t) \le u'(i)(t)$, $v(i)(t) \ge v'(i)(t)$, all $t$.
\end{corollary}
\textbf{Proof}
The result follows by applying Corollary \ref{C1} successively over $[0, \ta], [\ta, 2\ta], \cdots$ \hfill $\Box$

\begin{corollary} \label{C4}
With the same notation and hypothesis as in Corollary \ref{C3}, suppose $e(i)(t) = e'(i)(t)$, and define $z(i)(t) = x'(i)(t) - x(i)(t)$.  Then
\begin{eqnarray*}
\dot{x}(i)(t)&=& e(i)(t) - b(i)(t) + \sum b(j)(t - \ta (j,i))r(j,i) ,\\
\dot{x}'(i)(t)&=& e'(i)(t) - b'(i)(t) + \sum b'(j)(t - \ta (j,i))r(j,i) ,\\
\dot{z}(i)(t)&=& -[b'(i)(t) - b(i)(t)] + \sum [b'(j)(t- \ta (j,i)) - b(j)(t-\ta (j,i))] r(j,i) .
\end{eqnarray*}
Then $z(i)(t) \ge 0$ for all $t$.
\end{corollary}

\section{Periodic solution}\label{sec-periodic}
We now consider FT control.
Suppose  that the external arrivals $e(i)(t)$ and saturation flow rates $c(i)(t)$ are all periodic with the same period $T$.  Let
\begin{equation}
\eib= \frac{1}{T} \int_0^T e(i)(t) dt,\; \cib = \frac{1}{T} \int_0^T c(i)(t) dt.\label{20}
\end{equation}
We establish a necessary  condition for the existence of a periodic solution to \eqref{12}-\eqref{14} with period $T$.  Let
\begin{equation}
\bib = \frac{1}{T} \int_0^T b(i)(t) dt = \frac{1}{T} \int_0^T b(i)(t-\ta (i,j)) dt, \;  \yib = \frac{1}{T} \int_0^T y(i)(t) dt.\label{21}
\end{equation}
Using this notation in \eqref{14}, and integrating over $[0,T]$ for a periodic solution $x$ gives 
\begin{eqnarray*}
0&=& x(i) (T) - x(i)(0) = \eib- \cib + \sum_j \bjb r(j,i) + \yib \\
&=& \eib - \cib + \sum_j \cjb r(j,i) +\yib - \sum_j \yjb r(j,i),
\end{eqnarray*}
or, in vector form, denoting the routing matrix $R = \{r(i,j)\}$,
\[ 0 = \eb - [I-R^T] \cb + [I-R^T]\yb.\]
Since every vehicle eventually leaves, $[I-R^T]^{-1} = I + R^T + R^{2T} + \cdots \ge 0$ exists and so the preceding condition becomes
\begin{equation}
0 = [I-R^T]^{-1} \eb - \cb + \yb ,\label{22}
\end{equation}
so, for a periodic solution to exist, one must have
\[\cb = [I-R^T]^{-1} \eb + \yb \ge [I-R^T]^{-1} \eb.\]
We impose the  slightly stronger \textit{stability condition}: there exists $\epsilon> 0$ such that
\begin{equation}
 \cb > [I-R^T]^{-1} \eb + \epsilon \1, \label{23}
 \end{equation}
which says that on average the service rate for each queue exceeds the  total arrival rate.

\begin{remark} \label{R0} \
\rm  In the non-periodic case the stability condition \eqref{23} may be replaced by
\begin{equation}
\int_0^{NT} c(t) dt > [I-R^T]^{-1} \int_0^{NT} e(t) dt + \epsilon NT \1 , \mbox{  for all } N, \mbox{ for some } \epsilon > 0.\label{23A}
\end{equation}
Lemmas \ref{L1} and \ref{L2} hold under this stability condition.
\end{remark}

\begin{lemma}\label{L1}
Every trajectory $x(t)$ of vehicle queue lengths is bounded.
\end{lemma}
\textbf{Proof} \  Let $x, u, v, y$ be a solution of \eqref{15}-\eqref{18}.  Then
\begin{equation}
\dot{x}(t) = e(t) - c(t) + R^T b(t) + y(t) + \delta (t) ,\label{24}
\end{equation}
in which
\[\delta (i)(t) = \sum_j [b(j)(t-\ta (j,i)) - b(j)(t)] r(j,i).\]
So
\[ \int_s^t \delta (i) (r) dr = \sum_j r(j,i) \big[ \int_{s-\ta (j,i)}^s b(j)(s) - \int_{t - \ta (j,i)}^t b(j)(s) \big ] ds.\]
Since $\ta (j,i) \le \tb$, and $b(j)(s) \le c(j)(s)$ is bounded, it follows that $ |\int_s^t \delta (r) dr | $ is bounded for all $s, t$.  So
$ |\int_s^t \delta (r) dr |  \le d\1$ for some constant $d$.

We show that if $x(i) (t_0) > NT \bar{c }(i)$, then $x (i)(t_0+ NT) -x(i)(t_0) < 0$ for some constant $N > \frac{d}{T\epsilon}$.  This is sufficient to show that $x(i)(t)$ is bounded, since the queue-length change per period is bounded.  Suppose that $x(i)(t_0) > NT \bar{c }(i)$.  Since
$x(i)(t+T) - x(i)(t) \ge -T \bar{c }(i)$, so $x(i)(t) > 0$ for $t_0 \le t \le t_0 + NT$ and $y(i)(t) = 0$ for $t_0  \le t \le t_0 + NT$.  Integrating
both sides of \eqref{24} gives
\begin{eqnarray}
x(i)(t+NT) - x(i)(t) &\le& NT \bar{e}{(i)} + NT \sum_j r(j,i) \bar{c }(j) -NT \bar{c }(i) + d \label{26} \\
&\le& -NT \epsilon + d < 0 \label{27},
\end{eqnarray}
in which \eqref{26} uses the fact that $\bar{b}(j) \le \bar{c}(j)$ and \eqref{27} follows from the stability condition and the inequality
$N > \frac{d}{T\epsilon}$. \hfill $\Box$
\subsection{Effect of initial conditions}
Suppose the stability condition \eqref{23} holds.  By Lemma \ref{L1} the queue lengths are bounded,  so from \eqref{26}
\[\int _0^t y(r)dr \rightarrow \infty,\]
component wise.  Thus the cumulative unused service at every queue is unbounded.  This implies that, independent of the service discipline (whether first in first out or something else), so long as the discipline is work conserving (i.e., a queue is served if it is non-empty, see \eqref{2}), every vehicle in
the initial condition will eventually leave the network.

The \textit{state} of the network at time $t$ is the pair $(x(t), \beta(t))$ where $\beta (t)$ is the history of departures over time $ [t -\tb, t] $, i.e., $\beta(t)$ is the function: $s \in [t -\tb, t] \mapsto b(s)$.  We want to show that trajectories starting from two different initial conditions, say  $(x(0), \beta (0))$ and $(x'(0), \beta '(0))$, will eventually converge.  Because of the monotonicity property, Corollary \ref{C3}, we may take one of the initial conditions to be zero.
\begin{lemma} \label{L2}
Let $x(t), \beta (t)$ be the trajectory starting from $(x(0), \beta (0))$, and let $z(t), \beta' (t)$ be the trajectory starting from $(0,0)$. Then
the two trajectories converge:
\begin{equation}
\lim_{t \rightarrow\infty} |x(t) - z(t)| \rightarrow 0, \mbox{ and } \lim_{t \rightarrow \infty }\sup_ {s \in [t -\tb, t]}  [\int_0^s [b(r) - b'(r)] dr \rightarrow 0,
\label{28}
\end{equation}
\end{lemma}
\textbf{Proof} \ The proof relies on an intuitive argument.  Consider the trajectory $(x(t), \beta(t))$. Color the vehicles in the initial state $(x(0), \beta (0))$ red, and color all vehicles entering the network after time 0, black.   In each queue there will be black and red vehicles.  Change the service discipline in each queue so that all black vehicles are served ahead of every red vehicle.  Then the red vehicles do not interfere with the movement of black vehicles and so the number of black vehicles in the queues and along the links will be identical to $(z(t), \beta' (t))$.  On the other hand the total number  (red plus black)
vehicles in the queues and along the links will be identical to $(x(t), \beta (t))$.  Because of the stability condition, every vehicle in the initial queue
$x(0)$ will eventually leave the network, that is  $x(t) - z(t) \rightarrow 0$, as $t \rightarrow \infty$.  But then the second part of \eqref{28} follows. \hfill $\Box$

\begin{remark} \rm
Suppose the arrivals and service processes, $e$ and $c$, are stochastic  and the stability condition \eqref{23A} holds almost
surely.  Then the effect of the initial state on the queue length process  will disappear over time since \eqref{28} will hold almost surely along every sample path.
\end{remark}
\subsection{Existence of periodic solutions}
We prove the existence of a unique periodic solution to which all trajectories converge.
\begin{theo} \label{T4} There exists a unique periodic state trajectory $(x^*, \beta^* )$, with period $T$, to which every trajectory converges.
\end{theo}
\textbf{Proof}  \ Consider the Poincare map $F: (x(0), \beta(0)) \mapsto (x(T), \beta (T))$ and the iterates $(x(nT), \beta (nT)) = F(x((n-1)T, \beta ((n-1)T))$.
Take $x(0)=0, \beta(0)=0$.  Then by monotonicity $(x(T) \ge x(0) =0, \beta (T) \ge \beta (0) = 0)$, and hence by repeatedly using monotonicity we
get:
\[x(0) \le x(T) \le x(2T) \le \cdots, \quad \beta(0 ) \le \beta (T) \le \beta (2T) \le \cdots .\]
Thus this sequence of states is increasing.  By Lemma \ref{L1} the sequence is bounded, so it converges to the state say $(x^*, \beta^*)$.  By
Theorem \ref{T1}, $F$ is continuous, so $F((x^*, \beta^*)) = (x^*, \beta^*)$, and the trajectory from this state is periodic with period $T$. \hfill $\Box$
\begin{corollary} \label{C5}
In the periodic trajectory $x^*$ every queue clears in each period, i.e., for each $i$ there exists $t_i \in [0,T]$ such that $x(i)(t_i) = 0$.  In
each period, the cumulative unused service is $\yb = \cb -[I-R^T]^{-1} \eb$.
\end{corollary}
\textbf{Proof}  \ If $x(i)(t) > 0$ for $t \in [0,T]$, $y(i) (t) = 0$ for $t \in [0,T]$ and so $\yb (i) = 0$, which contradicts the stability condition
\eqref{22}, \eqref{23}. \hfill $\Box$
\subsection{Finite time convergence}
\begin{theo}\label{T5}
Suppose every vehicle leaves after visiting at most $K$ queues.  Then every trajectory converges to the periodic trajectory in finite time.
\end{theo}
\textbf{Proof}  \  Revisit the proof of Lemma \ref{L2}.  Consider the trajectory $(x(t), \beta (t))$ starting in the initial state $(x(0), \beta (0))$.  Color the vehicles in the initial state red, and the vehicles arriving after time 0, black.  Color all vehicles in the state starting in state $(0,0)$, black.  The red vehicles will potentially be served infinitely often in each queue and so they will all be gone after a finite time.  At that time the
two trajectories wil coincide and will agree with the periodic trajectory.  \hfill $\Box$

\begin{figure}[h!]
\centering
\includegraphics[width=5in]{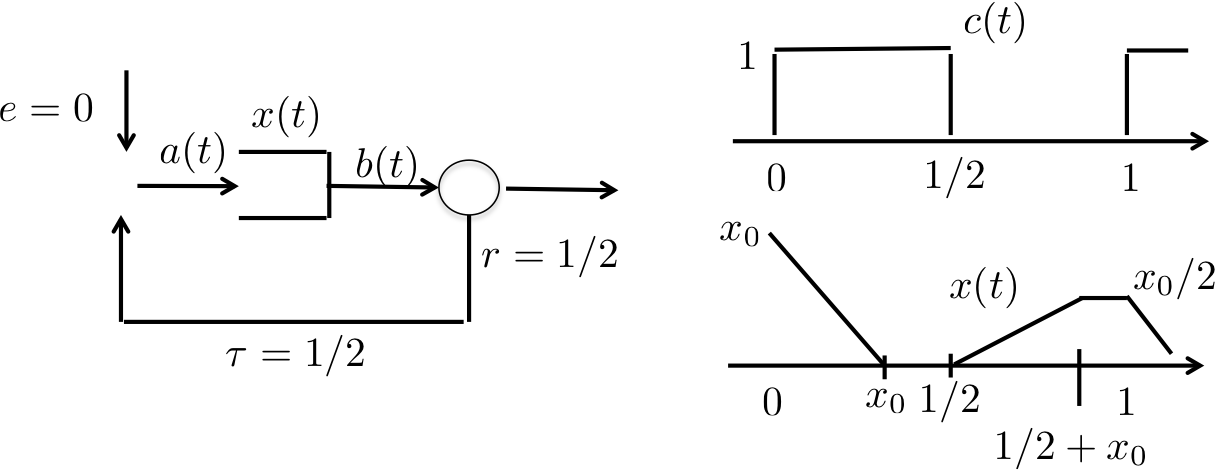}
\caption{
Vehicles recirculate (left); periodic service rate (right, top); queue $x(t)$ (right, bottom).}
\label{fig-example}
\end{figure}

\textbf{Example 2} \;
If vehicles can circulate indefinitely, convergence may take infinite time.  Figure \ref{fig-example} shows a network with a single queue with
initial size $x(0) = x_0$ and periodic service rate $c(t)$ with period 1, $c(t) = 1$ for $0\le t \le 1/2$ and $c(t) = 0$ for $1/2 < t \le 1$.
One-half ($r = 1/2$) of the departing vehicles return for service after travel time $\tau =1/2$; the other vehicles leave.
Suppose $x_0 < 1/2$.  Then all vehicles will depart by time $x_0$, one-half of them will leave and one-half or $x_0/2$ will re-enter the queue during
time $[1/2, 1/2 + x_0]$.  Since the service rate is 0 until time 1, so $x(t) = x_0/2$, $1/2 + x_0 < t <  1$.  At time 1, the queue is $x_0/2$
and there is no vehicle traveling in the link.   By induction, we have $x(n) = (1/2)^n x_0$, so convergence takes infinite time.

\begin{figure}[h!]
\centering
\includegraphics[width=5in]{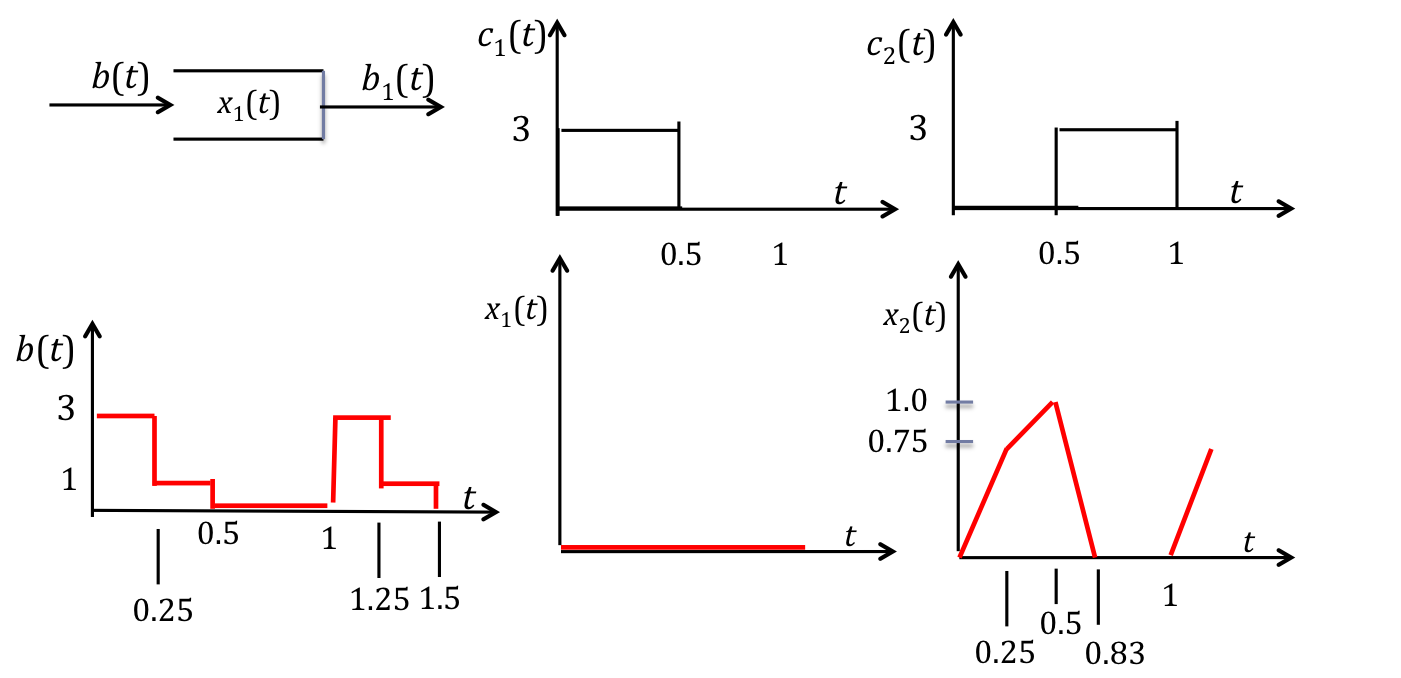}
\caption{
Departures $b(t)$ of Figure \ref{fig1} enter an intersection with service rate $c_1(t)$ or $c_2(t)$.}
\label{fig-example3}
\end{figure}

\textbf{Example 3} \;
From the periodic trajectory one can calculate performance measures such as average delay.  In Figure \ref{fig-example3} the departure process
$b(t)$ of Figure \ref{fig1} is the arrival process at the next intersection with the service rate $c_1(t)$ and queue $x_1(t) \equiv 0$
or service rate $c_2(t)$ and queue $x_2(t)$.  
The average delay per vehicle at the second intersection therefore  is
\[\int_0^1 x_1(t) dt = 0, \mbox{  or } \int_0^1 x_2(t) = \frac{23}{48} \approx 0.48.\]
Hence, depending on the offset of  the second signal, the average delay at the second intersection can take any value in [0, 0.48].  The average delay in the
first intersection with constant arrival can be calculated from the plot of $x(t)$ in Figure \ref{fig1} as
\[\int_0^1 x(t) dt = \frac{3}{16} \approx 0.19.\]
Webster's formula, commonly used as an approximation, gives the per vehicle delay as
\[d_w = \frac{T(1-g/T)^2}{2[1-(g/T)x]} + \frac{x^2}{2q(1-x)} - 0.65 \bigl (\frac{T}{q^2} \bigr)^{1/3} x ^{2 +5(g/T)} ,\]
with $T = 1$ is the cycle length, $g = 1/2$ is the green time, $q = \mbox{ flow} = 1$, $x = \mbox{ flow to capacity ratio} = 1/(g \times 3) = 2/3$, which works out as 
\[ d_w = \frac{3}{16} +\frac{2}{3} - 0.1 \approx 0.75.\]

\subsection{Finite storage capacity}
In the discussion so far it has been assumed that every queue has infinite storage capacity so service is never blocked.  We now modify this assumption.  Suppose queue $i$ has storage capacity $\xi (i)$.  This means that if $x(i)(t)$ reaches $\xi(i)$, additional vehicles cannot be accommodated and  service to queues upstream of queue $i$ is blocked.  More precisely, the service rate for queue $i$ is changed from
$c(i)(t)$ to $c(i)(t) s(i)(t)$ in which
\begin{equation}
 s(i)(t) = \left \{
\begin{array}{ll}
1, & \mbox{ if } x(j)(t) < \xi (j) \mbox{ for all } j \mbox{ such that } r(i,j) > 0 ,\\
0, & \mbox{ if } x(j) (t) = \xi (j) \mbox{ for some } j \mbox{ such that } r(i,j) > 0.
\end{array}
\right . \label{fifo}
\end{equation}
The system equations \eqref{14}-\eqref{18} remain the same except that $c(i)(t)$ is everywhere replaced by $c(i)(t) \times s(i)(t)$.
The formulation \eqref{fifo} implies that service from $i$ to $k$ is blocked even if $x(k)(t) < \xi (k)$.  This is a kind of first in first out
assumption.  We can avoid this by positing a separate queue for each movement as in \cite{MPtrc}.  The example below is unaffected in either case.

The possibility of blocking can destroy the previous results.  We can see this by examining again the single queue system of Figure \ref{fig-example}.  Suppose that  the queue has a storage capacity $\xi = 1$.  In this case the effective service rate is $c(t) \times s(t)$ and
$s(t) = \1 (x(t) < 1)$.  Suppose $e(t) \equiv 0$ and $x(t_0) = 1$. Then $x(t) \equiv 1$ and $b(t) \equiv 0$,  for $t \ge t_0$, and the system is in gridlock.  If $x(t_0) = 0$, then $x(t) \equiv 0$, for $t \ge t_0$ is another solution.

Suppose there is a constant arrival $e(t) \equiv \eb < 1/4$.  If $\xi = \infty$, the stability condition \eqref{23} holds and there will be
a periodic solution, $x^*$.   Suppose $\max_t x^*(t) = \bar{\xi}$.  If the storage capacity $\xi > \bar{\xi}$, then $x^*$ is also a solution.

Above, the service rate $c(i)(t) \times s(i)(t)$ is state-dependent and so the results above do \textit{not} apply.  In actuated traffic control, as opposed to FT control, the service rate indeed depends on the traffic state, so studying \eqref{14} for state-dependent service, $c(i, x(t), t)$, is important. The fundamental results on existence and monotonicity in case that $c(i, x, t)$ is Lipschitz in $x$ are obtained in \cite{rama-reflect}.
\section{Conclusion} \label{sec-conc}
A network of signalized intersections is modeled as a queuing network, whose state is the vector of queue lengths, together with the position of vehicles traveling along the links between intersections.  The state of the network evolves according to a delay-differential equation.
In this study each intersection is controlled by a fixed time controller with the same period or cycle.  External inputs are periodic with the same period.   Vehicles make turns at intersections in  fixed proportions.  The network is undersaturated.
The major conclusion of the study is that there is a unique periodic trajectory, which is globally asymptotically stable, that is, every trajectory converges to the periodic trajectory.  In case every vehicle leaves the network after traveling over a bounded number of links, the convergence occurs in finite time.

From the periodic trajectory one can easily calculate every possible performance measure such as delay, travel time, amount of time
service is wasted, and progression quality.  Thus an important question for future research is to find an algorithm to calculate the
periodic orbit.  Another question is to study the behavior of traffic in networks with actuated control in which the service rate is state-dependent.  It seems a reasonable conjecture that if the control is a function of the queue lengths (as in \cite{MPtrc}), there will again
be a unique asymptotically stable trajectory.

\bibliographystyle{plainnat}

\bibliography{ft_additional_bib,/Users/varaiya/Documents/varaiya-Main/varaiya/Bib/traffic,/Users/varaiya/Documents/varaiya-Main/Papers.dir/Bib/traffic}

\end{document}